\documentclass[psamsfonts,12pt]{amsart}
\usepackage{verbatim}
\usepackage{amsmath}
\usepackage{amssymb}
\usepackage{epsfig}

\theoremstyle{plain}
   \newtheorem{theorem}{Theorem}[section]
\newtheorem*{theorem*}{The Livshitz Theorem}

   \newtheorem{lemma}[theorem]{Lemma}
   \newtheorem{corollary}[theorem]{Corollary}
   \newtheorem{prop}[theorem]{Proposition}

   \theoremstyle{definition}
   \newtheorem*{definition}{Definition}

   \theoremstyle{definition}
   \newtheorem*{definitions}{Definitions}

   \theoremstyle{remark}
   \newtheorem*{remark}{Remark}

\theoremstyle{remark}

\newcommand{\R}{{\mathbb R}}

\makeatletter
   \def\@setcopyright{}
   \def\serieslogo@{}
   \makeatother

\newcommand {\C}{{\mathbb C}}

\newcommand {\T}{{\mathbb T}}
\newcommand {\Z}{{\mathbb Z}}
\newcommand {\Q}{{\mathbb Q}}

\newcommand {\Rd} {{\mathbb R}^d}

\newcommand {\Zd} {{\mathbb Z} ^d}
\newcommand{\bn} {{\mathbf n}}

\newcommand {\g}{{\gamma}}
\newcommand {\G}{{\Gamma}}
\newcommand {\lm}{{\lambda}}

\newcommand {\al}{{\alpha}}
\newcommand {\la} {{\lambda}}

\newcommand {\GnZ} {{GL(n,{\mathbb Z})}}
\newcommand {\GnR} {{GL(n,{\mathbb R})}}
\newcommand {\GnQ} {{GL(n,{\mathbb Q})}}
\newcommand {\SnZ} {{SL(n,{\mathbb Z})}}
\newcommand {\SnR} {{SL(n,{\mathbb R})}}

\newcommand {\Id}{{\rm Id}}
\newcommand {\id}{{\rm id}}
\newcommand{\Fix}{{\rm Fix}}

\newcommand{\foot}[1]{\mbox{}\marginpar{\raggedleft\hspace{0pt}
\Tiny
#1}}

\begin{document}\frenchspacing
\author{Anatole Katok}
   \address{A. Katok: Department of Mathematics, Pennsylvania
State Univer\-sity,
   University Park, PA 16802, USA}
   \email{katok\_a@math.psu.edu}

\author{Svetlana Katok}
   \address{S. Katok: Department of Mathematics, Pennsylvania
State Univer\-sity,
   University Park, PA 16802, USA}
   \email{katok\_s@math.psu.edu}

\author{Klaus Schmidt}
\address{K. Schmidt: Mathematics Institute, University of
Vienna, Strudl\-hofgasse~4, A-1090 Vienna, Austria,
\newline\indent
\textup{and} \newline\indent Erwin Schr\"odinger Institute for
Mathematical Physics, Boltzmann\-gasse~9, A-1090 Vienna, Austria}
\email{klaus.schmidt@univie.ac.at}

\title[Measure-theoretic rigidity]{Rigidity of measurable
structure for $\Z^d$--actions by automorphisms of a torus}

\date{\today}
\begin{abstract}We show that
for certain  classes of actions of $\Z^d,\,\,d\ge 2$, by 
automorphisms of the torus any measurable conjugacy has to be affine, 
hence measurable conjugacy implies algebraic conjugacy; similarly any 
measurable factor is algebraic, and algebraic and affine centralizers 
provide invariants of measurable conjugacy. Using the algebraic 
machinery of dual modules and information about class numbers of 
algebraic number fields we consruct various examples of 
$\Z^d$-actions by Bernoulli automorphisms  whose measurable orbit 
structure is rigid, including actions which are weakly isomorphic but 
not isomorphic.  We show that the  structure of the centralizer for 
these actions may or may not serve as a distinguishing 
measure--theoretic invariant.

\end{abstract} \subjclass{}
\keywords{}
\thanks{The research of the first author was partially
supported by NSF grant DMS-9704776. The first two authors are
grateful to the Erwin Schr\"odinger Institute, Vienna, and the
third author to the Center for Dynamical Systems at Penn State
University, for hospitality and support while some of this
work was done.}

   \maketitle
\section {Introduction; description of results}

In the course of the last decade various rigidity properties
have been found for two different classes of actions by
higher--rank abelian groups: on the one hand, certain Anosov and
partially hyperbolic actions of $\Z^d$ and $\R^d,\,\, d\ge 2$,
on compact manifolds (\cite{KS1, KS2, KS4}) and,
on the other, actions of
$\Z^d,\,\, d\ge 2$, by automorphisms of compact abelian groups
(cf. e.g. \cite{KSch, KiS}).
Among these rigidity phenomena is a relative scarcity
of invariant measures which stands in contrast with the
classical case $d=1$ (\cite{KS3}).

In this paper we make the first
step in investigating a different albeit related phenomenon:
rigidity of the measurable orbit structure with respect to
the natural smooth invariant measure.

In the classical case of actions by
$\Z$ or $\R$ there are certain natural classes of
measure--preserving transformations which possess such rigidity:
ergodic translations on compact abelian groups give a rather
trivial example, while horocycle flows and other homogeneous
unipotent systems present a much more interesting one
\cite{R1, R2, R3}. In contrast to such situations,
individual elements of the higher--rank actions mentioned above
are Bernoulli automorphisms. The measurable orbit structure of a
Bernoulli map can be viewed as very ``soft''. Recall that the
only metric invariant of Bernoulli automorphisms is entropy
(\cite{O}); in particular, weak
isomorphism is equivalent to isomorphism for Bernoulli maps
since it implies
equality of entropies. Furthermore, description of
centralizers, factors, joinings and other invariant objects
associated with a Bernoulli map is impossible in reasonable
terms since each of these objects is huge and does not possess
any discernible structure.

In this paper we demonstrate that some very natural actions
of $\Z^d,\,\,d\ge 2$, by
Bernoulli automorphisms display a remarkable
rigidity of their measurable orbit structure. In particular,
isomorphisms between such actions, centralizers, and factor maps
are very restricted, and a lot of algebraic
information is encoded in the measurable
structure of such actions (see Section \ref{s:rigidity}).

All these properties occur for broad subclasses of both main
classes of actions of higher--rank abelian groups mentioned
above: Anosov and partially hyperbolic actions on compact
manifolds, and actions by automorphisms of compact abelian
groups.
   However, at present we are unable to present sufficiently
definitive general results due to various
difficulties of both conceptual and technical nature. Trying to
present the most general available results would lead to
cumbersome notations and inelegant formulations. To avoid that
we chose to restrict our present analysis to a smaller class
which in fact represents the intersection of the two, namely
the actions of
$\Z^d,\,\,d\ge 2$, by automorphisms of the torus. Thus we study
the measurable structure of such actions with respect to
Lebesgue (Haar) measure from the point of view of ergodic
theory.

Our main purpose is to demonstrate several striking phenomena
by means of applying to specific examples general rigidity
results  which are presented in Section \ref{s:rigidity} and are
based on rigidity of invariant measures developed in
\cite{KS3} (see \cite{KaK} for further results along these
lines including rigidity of joinings). Hence we do not strive
for the greatest possible generality even within the class of
actions by automorphisms of a torus. The basic algebraic setup
for irreducible actions by automorphisms of a torus is
presented in Section \ref{s:irred}.
Then we adapt further necessary
algebraic preliminaries to the special but in a sense most
representative case of Cartan actions, i.e. to
$\Z^{n-1}$--actions by hyperbolic automorphisms of the
$n$--dimensional torus (see
Section \ref{s:ca}).

The role of entropy for a smooth action of a higher--rank abelian
group $G$ on a finite-dimensional manifold is played by the
{\it entropy function} on $G$ whose
values are entropies of individual elements of the action (see
Section \ref{ss:high-rank} for more details) which is
naturally invariant of isomorphism and also of weak isomorphism
and is equivariant with respect to a time change.

In Section \ref{s:examples}
we produce several kinds of specific examples of actions by
ergodic (and hence Bernoulli) automorphisms of tori with
the same entropy function. These examples provide concrete
instances when general criteria developed in Section
\ref{s:rigidity} can be applied. Our examples include:
\begin{itemize}
\item[(i)]
actions which are not weakly isomorphic
(Section \ref{ss:not weakly}),
\item[(ii)]
actions which are weakly isomorphic but not isomorphic,
such that one action is a maximal action by Bernoulli
automorphisms and the other is not (Section \ref{ss:cent}),
\item[(iii)] weakly
isomorphic, but nonisomorphic, maximal actions
   (Section \ref{ss:max-Cart}).
\end{itemize}

Once rigidity of conjugacies is established, examples of type
(i) appear in a rather simple--minded fashion: one simply
constructs actions with the same entropy data which are not
isomorphic over $\Q$. This is not surprising since entropy
contains only partial information about eigenvalues. Thus one
can produce actions with different eigenvalue structure but
identical entropy data.

Examples of weakly isomorphic but nonisomorphic actions
are more sophisticated. We find them among
Cartan actions (see Section \ref{s:ca}).
The centralizer of a Cartan action in the group of
automorphisms of the torus is (isomorphic to) a finite extension
of the acting group, and in
some cases Cartan actions isomorphic over $\Q$ may be
distinguished by looking at the index of the group in its
centralizer (type (ii); see Examples 2a and 2b). The underlying
cause for this phenomenon is the existence of algebraic number fields
$K=\Q(\lambda)$, where
$\lambda$ is a unit, such that the ring of integers $\mathcal
O_K\ne\Z[\lambda]$. In general finding even simplest possible
examples for $n=3$ involves the use of data from algebraic
number theory and rather involved calculations.
   For examples of type
(ii) one may use some special tricks which allow to find some
of these and to show nonisomorphism without a serious
use of symbolic manipulations on a computer.

A Cartan action $\alpha$ of $\Z^{n-1}$ on $\T^n$ is called
maximal if its centralizer in the group of automorphisms of
the torus is equal to
$\alpha(\Z^{n-1})\times\{\pm\Id\}$. A maximal Cartan action
turns out to me maximal in the above sense: it cannot be
extended to any action of a bigger abelian group
by Bernoulli automorphisms.

Examples of maximal Cartan actions isomorphic over
$\Q$ but not isomorphic (type (iii)) are the most remarkable.
Conjugacy over
$\Q$ guarantees that the actions by automorphisms of the torus
$\T^n$ arising from their centralizers are weakly isomorphic
with finite fibres. The mechanism providing
obstructions for algebraic isomorphism in this case involves the
connection between the class number of an algebraic number
field and $GL(n,\Z)$--conjugacy classes of matrices in
$SL(n,\Z)$ which have the same characteristic polynomial (see
Example 3). In finding these examples the use of
computational number--theoretic algorithms (which in our
case were implemented via the Pari-GP package) has been essential.

   One of our central conclusions is that for a broad
class of actions of $\Z^d,\,\,d\ge 2$, (see condition
$(\mathcal R)$ in Section \ref{ss:high-rank}) the conjugacy class
of the centralizer of the action in the group of {\it affine}
automorphisms of the torus is an invariant of measurable
conjugacy. Let $Z_\mathit{meas}(\al)$ be the centralizer of the
action
$\al$ in the group of measurable automorphisms. As it turns out
in all our examples but Example 3b, the conjugacy class of the pair
$(Z_\mathit{meas}(\al),\al)$ is a distinguishing invariant of the
measurable isomorphism.  Thus, in particular, Example 3b shows that
there are weakly isomorphic, but nonisomorphic actions
for which the affine and hence the measurable centralizers are
isomorphic as abstract groups.

We would like to acknowledge a contribution of J.-P.
Thouvenot to the early development of ideas which led to
this paper. He made an important observation that rigidity of
invariant measures can be used to prove rigidity of isomorphisms
via a joining construction (see Section \ref{ss:conj}).

\section {Preliminaries}

\subsection {Basic ergodic theory} Any invertible (over $\Q$)
integral $n\times n$ matrix
$A\in M(n,\Z)\cap GL(n,\Q)$ determines an endomorphism of the
torus
$\T^n=\R^n/\Z^n$ which we denote by $F_A$. Conversely, any
endomorphism of $\T^n$ is given by a matrix from
$A\in M(n,\Z)\cap GL(n,\Q)$. If, in addition, $\det A=\pm
1$, i.e. if $A$ is invertible over $\Z$, then $F_A$ is an
automorphism of $\T^n$ (the group of all such
$A$ is denoted by $GL(n,\Z)$). The map $F_A$ preserves
Lebesgue (Haar) measure
$\mu$; it is
ergodic with respect to
$\mu$ if and only if there are no roots of unity among the
eigenvalues of $A$, as was first pointed out by Halmos
(\cite{H}). Furthermore, in this case there are eigenvalues of
absolute value greater than one and
$(F_A,\la)$ is an exact endomorphism. If $F_A$ is an automorphism
it is in fact Bernoulli (\cite{Kat}). For
simplicity we will call such a map
$F_A$ an {\it ergodic toral endomorphism} (respectively, {\it
automorphism}, if $A$ is invertible). If all eigenvalues of $A$
have absolute values different from one we will call the
endomorphism (automorphism) $F_A$ {\it hyperbolic}.

   When it does not lead to
a confusion we will not distinguish between a matrix $A$ and
corresponding toral endomorphism $F_A$.

   Let $\lm_1,\dots , \lm_n$ be the eigenvalues of the matrix
$A$, listed with their multiplicities. The entropy
$h_{\mu}(F_A)$ of $F_A$ with respect to Lebesgue measure is
equal to
$$\sum_{\{i:|\lm_i|>1\}}\log|\lm_i|.$$ In particular, entropy
is determined by the conjugacy class of the matrix $A$ over $\Q$ (or
over
$\C$). Hence {\it all ergodic toral automorphisms which are conjugate
over $\Q$ are measurably conjugate with respect to Lebesgue
measure.}

Classification, up to a conjugacy over $\Z$, of matrices in
$\SnZ$, which are irreducible and conjugate over $\Q$ is
closely related to the notion of class number of an algebraic
number field. A detailed
discussion relevant to our purposes appears in Section
\ref{ss:LM}. Here we only mention the simplest case
$n=2$ which is not directly related to rigidity. In
this case trace determines conjugacy class over $\Q$ and, in
particular, entropy. However if the class number of the
corresponding number field is greater than one there are
matrices with the given trace
which are not conjugate over
$\Z$. This algebraic distinctiveness is not reflected in the
measurable structure: in fact, in the case of equal entropies
the classical Adler--Weiss
construction of the Markov partition in \cite{AW} yields metric
isomorphisms which are more concrete and specific than in the
general Ornstein isomorphism theory and yet not algebraic.

\subsection{Higher rank actions}\label{ss:high-rank}
Let $\al$ be an action by
commuting toral automorphisms given by integral matrices
$A_1,\dots,A_d$. It defines an embedding
$\rho_\al:\Z^d\to GL(n,\Z)$ by
\[
\rho_\al^\bn=A_1^{n_1}\dots A_d^{n_d},
\] where $\bn=(n_1,\dots,n_d)\in \Z^d$, and we have
\[
\al^\bn=F_{\rho_{\al}^\bn}.
\] Similarly, we write  $\rho_{\al}:\Z^d_+\to M(n,\Z)\cap
GL(n,\Q)$ for an action by endomorphisms. Conversely, any
embedding $\rho:\Z^d\to GL(n,\Z)$ (respectively,
$\rho:\Z^d_+\to M(n,\Z)\cap GL(n,\Q))$ defines an action by
automorphisms (respectively, endomorphisms) of $\mathbb{T}^n$
denoted by $\al_{\rho}$.

Sometimes we will not explicitly distinguish between an action
and the corresponding embedding, e.g. we may talk about ``the
centralizer of an action in $GL(n,\Z)$'' etc.

\begin{definitions}

Let $\alpha$ and
$\alpha'$ be two actions of $\Z^d$ ($\Z^d_+$) by automorphisms
(endomorphisms) of $\T^n$ and $\T^{n'}$, respectively.
The actions $\alpha$ and $\alpha'$ are {\em measurably} (or
{\em metrically}, or {\em measure--theoretically})
{\em isomorphic} (or {\em conjugate})
if there exists a Lebesgue measure--preserving bijection
$\varphi:\T^n\to\T^{n'}$ such that
$\varphi\circ\alpha=\alpha'\circ\varphi$.

The actions $\alpha$ and $\alpha'$ are {\em measurably isomorphic up
to a time change} if there exist a measure--preserving
bijection
$\varphi:\T^n\to\T^{n'}$ and a $C\in GL(d,\Z)$ such that
$\varphi\circ\alpha\circ C=\alpha'\circ\varphi$.

The action $\alpha'$ is a {\em measurable factor} of $\alpha$ if
there exists a Lebesgue measure--preserving transformation
$\varphi:\T^n\to\T^{n'}$ such that
$\varphi\circ\alpha=\alpha'\circ\varphi$. If, in particular,
$\varphi$ is almost everywhere finite--to--one, then $\alpha'$
is called a {\em finite factor} or a {\em factor with finite
fibres} of $\alpha$.

Actions $\alpha$ and $\alpha'$ are {\em weakly
measurably isomorphic} if each is a measurable factor of the other.

A {\em joining} between $\alpha$ and $\alpha'$ is a
measure $\mu$ on $\T^n\times\T^{n'}=\T^{n+n'}$ invariant under
the Cartesian product action $\alpha\times\alpha'$ such that
its projections into $\T^n$ and $\T^{n'}$ are Lebesgue
measures. As will be explained in Section \ref{s:rigidity},
conjugacies and factors produce special kinds of joinings.
\end{definitions}

These measure--theoretic notions have natural algebraic
counterparts.

\begin{definitions}
The actions $\alpha$ and $\alpha'$
are {\em algebraically
isomorphic} (or {\em conjugate}) if $n=n'$ and if
there exists a group automorphism
$\varphi:\T^n\to\T^{n}$ such that
$\varphi\circ\alpha=\alpha'\circ\varphi$.

The actions $\alpha$ and $\alpha'$ are {\em algebraically isomorphic
up to a time change} if there exists an automorphism
$\varphi:\T^n\to\T^{n}$ and $C\in GL(d,\Z)$ such that
$\varphi\circ\alpha\circ C=\alpha'\circ\varphi$.

The action $\alpha'$ is an {\em algebraic factor} of $\alpha$
if there exists a surjective homomorphism
$\varphi:\T^n\to\T^{n'}$ such that
$\varphi\circ\alpha=\alpha'\circ\varphi$.

The actions $\alpha$ and $\alpha'$ are {\em weakly
algebraically isomorphic} if each is an algebraic factor of the other.
In this case $n=n'$ and each factor map has finite
fibres.

Finally, we call a map $\varphi: \T^n\to \T^{n'}$ {\em
affine} if there is a surjective continuous group homomorphism
$\psi: \T^n\to \T^{n'}$ and $x'\in\T^{n'}$ s.t.
$\varphi(x)=\psi(x)+x'$ for every $x\in \T^n$.
\end{definitions}

As already mentioned, we intend to show that under certain
condition for $d\ge 2$, measure theoretic properties imply
their algebraic counterparts.

We will say that an algebraic factor
$\al'$ of $\al$ is a {\it rank--one factor} if $\al'$ is an
algebraic factor of $\al$ and $\al'(\Z^d_+)$
contains a cyclic sub--semigroup of finite index.

The most general situation when
certain rigidity phenomena appear is the following :

\vspace {.1in}
$(\mathcal R')$: {\em The action $\al$ does not
possess nontrivial rank--one algebraic factors.}
\vspace {.1in}

In the case of actions by automorphisms the condition
$(\mathcal R')$ is equivalent to the following condition
$(\mathcal R)$ (cf. \cite{S}):

\vspace {.1in}
$(\mathcal R)$: {\em The action $\al$ contains a
group, isomorphic to $\Z^2$, which consists of ergodic
automorphisms. }

By Proposition 6.6 in \cite{Sch2}, Condition $(\mathcal{R})$ is
equivalent to saying that {\em the restriction of $\alpha $
to a subgroup isomorphic to $\Z^2$ is mixing.}
\vspace {.1in}

{\it A Lyapunov exponent} for an action $\al$ of $\Z^d$ is a
function $\chi: \Z^d\to \R$ which associates to each
${\bf n}\in\Z^d$
the logarithm of the absolute value of the eigenvalue for
$\rho_\al^\bn$ corresponding to a fixed eigenvector. Any
Lyapunov exponent is a linear function; hence it extends
uniquely to $\R^d$. The {\it multiplicity} of an exponent is
defined as the sum of multiplicities of eigenvalues corresponding to
  this exponent. Let
$\chi_i,\,\,i=1,\dots, k$, be the different Lyapunov exponents and let
$m_i$ be the multiplicity of $\chi_i$. Then the
entropy formula for a single toral endomorphism implies that
$$h_{\al}(\bn)=h_\mu(\rho_\al^\bn)=\sum_{\{i: \chi_i({\bf n})>
0\}}m_i\chi_i(\bf n).$$

The function $h_{\al}: \Z^d\to \R$ is called {\it the
entropy function} of the action $\al$. It naturally extends
to a symmetric, convex piecewise linear function of $\R^d$.
Any cone in $\R^d$ where all Lyapunov exponents have constant
sign is called a {\it Weyl chamber}. The entropy function is
linear in any Weyl chamber.

The entropy function is a prime invariant of measurable
isomorphism; since entropy does not increase for factors
the entropy function is also invariant of a weak measurable
isomorphism.
   Furthermore it changes
equivariantly with respect to automorphisms of $\Z^d$.

\begin{remark} it is interesting to point out that the
convex piecewise linear structure of the entropy function
persists in much greater generality, namely for smooth actions on
differentiable manifolds with a Borel invariant measure with
compact support.
\end{remark}

\subsection{Finite algebraic factors and invariant lattices}
Every algebraic action has many algebraic factors with finite
fibres. These factors are in one--to--one correspondence with
lattices $\G\subset\R^n$ which contain the standard lattice
$\G_0=\Z^n$, and which satisfy that
$\rho_{\al}(\G)\subset\G$.
The factor--action associated with a particular lattice
$\G\supset \G_0$ is denoted by
$\al_{\G}$. Let us point out that in the case of actions by
automorphisms such factors are also invertible:
if $\G\supset\G_0$ and $\rho_{\al}(\G)\subset\G$, then
$\rho_{\al}(\G)=\G$.

Let $\G\supset\G_0$ be a lattice. Take
any basis in $\G$ and let $S\in GL(n,\Q)$ be the matrix which
maps the standard basis in $\G_0$ to this basis. Then obviously
the factor--action $\al_{\G}$ is equal to the action
$\al_{S\rho_{\al}S^{-1}}$. In particular,
$\rho_{\al}$ and
$\rho_{\al_{\G}}$ are conjugate over $\Q$, although not
necessarily over $\Z$. Notice that conjugacy over $\Q$ is
equivalent to conjugacy over $\R$ or over $\C$.

For any positive integer $q$, the lattice $\frac
1{q}\G_0$ is invariant under any  automorphism in $GL(n,\Z)$
and gives rise to a factor which is conjugate to
the initial action: one can set
$S=\frac 1{q}\Id$ and obtains that $\rho_{\al}=\rho_{\al_{\frac
1{q}\G_0}}$. On the other hand one can find, for any lattice
$\G\supset\G_0$, a positive integer
$q$ such that $\frac 1{q}\G_0\supset\G$ (take $q$ the least
common multiple of denominators of coordinates for a basis of
$\G$). Thus $\al_{\frac 1{q}\G_0}$ appears as a factor of
$\al_{\G}$. Summarizing, we have the following properties
of finite factors.

\begin{prop}\label{prop-finfactors}Let $\al$ and $\al'$ be
$\Z^d$--actions by automorphism of the torus $\T^n$. The
following are equivalent.
\begin{enumerate}
\item $\rho_{\al}$ and $\rho_{\al'}$ are conjugate over $\Q$;
\item there exists an action $\al''$ such that both $\al$ and
$\al'$ are isomorphic to finite algebraic factors of $\al''$;
\item $\al$ and $\al'$ are weakly algebraically isomorphic,
i.e. each of them is isomorphic to a finite algebraic factor of
the other.
\end{enumerate}
\end{prop}

Obviously, weak algebraic isomorphism implies weak
measurable isomorphism. For
$\Z$--actions by Bernoulli automorphisms, weak isomorphism
implies isomorphism since it preserves entropy, the only
isomorphism invariant for Bernoulli maps. In Section
\ref{s:rigidity} we will show that, for  actions by toral
automorphisms satisfying Condition $(\mathcal R)$, measurable isomorphism
implies algebraic isomorphism. Hence, existence of such actions
which are conjugate over $\Q$ but not over $\Z$ provides
examples of actions by Bernoulli maps which are weakly
isomorphic but not isomorphic.

\subsection{Dual modules} For any action $\al$ of $\Z^d$ by
automorphisms of a compact abelian group $X$ we denote by $\hat\al$ the
dual action on the discrete group $\hat X$ of characters of $X$.
For an element $\chi\in \hat X$ we denote $\hat X_{\al_,\chi}$ the
subgroup of
$\hat X$ generated by the orbit $\hat \al\chi$.

\begin{definition} The
action $\al$ is called {\it cyclic} if $\hat X_{\al_,\chi}=\hat X$
for some $\chi\in\hat X$.
\end{definition}

Cyclicity is obviously an invariant of
algebraic conjugacy of actions up to a time change.

More generally, the dual group $\hat X$ has the structure of a
module over the ring $\Z[u_1^{\pm1},\dots,u_d^{\pm1}]$ of
Laurent polynomials in $d$ commuting variables. Action by the
generators of $\hat \al$ corresponds to multiplications by
independent variables. This module is called {\it the dual
module} of the action $\al$ (cf.\ \cite{Sch1, Sch2}). Cyclicity
of the action corresponds to the condition that this module
has a single generator. The structure of the
dual module up to isomorphism is an invariant of algebraic
conjugacy of the action up to a time change.

In the case of the torus $X=\T^n$ which concerns us in this
paper one can slightly modify the construction of the dual
module to make it more geometric. A $\Z^d$-action
$\al$ by automorphisms of the torus $\R^n/\Z^n$ naturally
extends to an action on $\R^n$ (this extension coincides with
the embedding $\rho_{\al}$ if matrices are
identified with linear transformations). This action preserves
the lattice $\Z^n$ and furnishes $\mathbb{Z}^n$ with the
structure of a module over the ring
$\Z[u_1^{\pm1},\dots,u_d^{\pm1}]$. This module is --- in
an obvious sense --- a \emph{transpose} of the
dual module defined above. In particular, the condition of
cyclicity of the action does not depend on which of these
two definitions of dual module one adopts.

\subsection{Algebraic and affine centralizers} Let
$\al$ be an action of $\Z^d$ by toral automorphisms, and let
$\rho_\al(\Z^d)=\{\rho_\al^\bn: n\in\Z^d\}$.
The \emph{centralizer} of $\al$ in the group of
automorphisms of $\T^n$ is denoted by
$Z(\al)$ and is not distinguished from the
centralizer of $\rho_\al(\Z^d)$ in $\GnZ$.

Similarly, the centralizer of $\al$ in the semigroup of all
endomorphisms of
$\T^n$ (identified with the centralizer of $\rho_\al(\Z^d)$ in
the semigroup $M(n,\Z)\cap\GnQ$) is
denoted by $C(\al)$.

The centralizer of $\al$ in the group of affine automorphisms
of $\T^n$ will be denoted by $Z_\mathit{Aff}(\al)$.

The centralizer of $\al$ in the semigroup of surjective affine
maps of
$\T^n$ will be denoted by $C_\mathit{Aff}(\al)$.

\section {Irreducible actions}\label{s:irred}

\subsection{Definition}

The action
$\al$ on $\mathbb{T}^n$ is called {\it irreducible} if any
nontrivial algebraic factor of $\al$ has finite fibres.

The following characterization of irreducible actions is
useful (cf. \cite{B}).

\begin{prop} \label{Berend}
The following
conditions are equivalent:
\begin{enumerate}
\item $\al$ is irreducible;
\item $\rho_{\al}$ contains a matrix with characteristic
polynomial irreducible over $\Q$;
\item $\rho_{\al}$ does not have a nontrivial invariant
rational subspace or, equivalently, any
$\al$--invariant closed subgroup of $\T^n$ is finite.
\end{enumerate}
\end{prop}

\begin{corollary}\label{cor:irred} Any irreducible action $\alpha$
of $\Z^d_+,\,\,d\ge 2$, satisfies condition \textup{($\mathcal R'$)}.
\end{corollary}
\begin{proof}
A rank one algebraic factor has to have fibres of positive
dimension. Hence the pre--image of the origin under the factor
map is a union of finitely many rational tori of positive
dimension and by Proposition \ref{Berend} $\alpha$ cannot be
irreducible.
\end{proof}

\subsection{Uniqueness of cyclic actions}
Cyclicity uniquely determines an irreducible action up
to algebraic conjugacy within a class of weakly algebraically
conjugate actions.

\begin{prop}\label{prop:cyclic}If $\al$ is an irreducible cyclic
action of
$\Z^d,\,\,d\ge 1$, on $\T^n$ and $\al'$ is another cyclic action
such that $\rho_\al$ and $\rho_{\al'}$ are conjugate over $\Q$,
then $\al$ and $\al'$ are algebraically isomorphic.
\end{prop}

For the proof of Proposition \ref{prop:cyclic} we need an elementary lemma.

\begin{lemma}Let $\rho:\Z^d\to \GnZ$ be an irreducible
embedding. The centralizer of $\rho$ in $\GnQ$ acts
transitively on $\Z^n\setminus\{0\}$.
\end{lemma}

\begin{proof} By diagonalizing $\rho$ over $\C$ and taking the
real form of it, one immediately sees that the centralizer of
$\rho$ in $\GnR$ acts transitively on vectors with nonzero
projections on all eigenspaces and thus has a single open and dense orbit.
Since the centralizer over
$\R$ is the closure of the centralizer over $\Q$, the $\Q$-linear
span of the orbit of any integer or rational vector under the
centralizer is an invariant rational subspace. Hence any
integer point other than the origin belongs to the single open
dense orbit of the centralizer of $\rho $ in $\GnR$.
This implies the statement of the lemma.
\end{proof}

\begin{proof}[Proof of Proposition \ref{prop:cyclic}] Choose $C\in M(n,\Z)$
such that $C\rho_{\al'}C^{-1}=\rho_\al$. Let $\mathbf k,\mathbf l\in \Z^n$ be
cyclic vectors for $\rho_\al|_{\Z^n}$ and $\rho_{\al'}|_{\Z^n}$,
respectively.

Now consider the integer vector $C(\mathbf l)$ and find $D\in
\GnQ$ commuting with $\rho_\al$ such that $DC(\mathbf l)=\mathbf
k$. We have $DC\rho_{\al'}C^{-1}D^{-1}=\rho_\al$. The
conjugacy $DC$ maps bijectively the $\Z$--span of the
$\rho_{\al'}$--orbit of $\mathbf l$ to $\Z$--span of the
$\rho_\al$--orbit of $\mathbf k$. By cyclicity both spans
coincide with $\Z^n$, and hence $DC\in \GnZ$.
\end{proof}

\subsection{Centralizers of integer matrices and algebraic
number fields}\label{ss:irr-integers}

There is an intimate connection between irreducible actions on
$\T^n$ and groups of units in number fields of
degree $n$. Since this connection (in the particular case
where the action is Cartan and hence the number field is
totally real) plays a central role in the construction of our
principal examples (type (ii) and (iii) of the Introduction),
we will describe it here in detail even though most of this
material is fairly routine from the point of view of
algebraic number theory.

Let $A\in\GnZ$ be a matrix with an irreducible characteristic
polynomial $f$ and hence distinct eigenvalues. The centralizer of
$A$ in
$M(n,\Q)$ can be identified with the ring of all polynomials in
$A$ with rational coefficients modulo the principal ideal
generated by the polynomial $f(A)$, and hence with the field
$K=\Q(\la)$, where
$\la$ is an eigenvalue of $A$, by the map
\begin{equation}\label{eq:gamma}
\g: p(A)\mapsto p(\la)
\end{equation}
with $p\in\Q[x]$.
Notice that if $B=p(A)$ is an integer matrix then $\g(B)$
is an algebraic integer, and if $B\in GL(n,\Z)$ then $\g(B)$
is an algebraic unit (converse is not
necessarily true).

\begin{lemma}\label{l:inj} The map $\g$ in \eqref{eq:gamma} is
injective.
\end{lemma}

\begin{proof} If $\g(p(A))=1$ for $p(A)\ne \Id$, then $p(A)$
has $1$ as an eigenvalue, and hence has a rational
subspace consisting of all invariant vectors. This subspace
must be invariant under $A$ which contradicts its
irreducibility.
\end{proof}

Denote by $\mathcal O_K$ the ring of integers in $K$, by
$\mathcal U_K$ the group of units in $\mathcal O_K$, by
$C(A)$ the centralizer of
$A$ in
$M(n,\Z)$ and by $Z(A)$ the centralizer of
$A$ in the group
$GL(n,\Z)$.

\begin{lemma}\label{l:rings-units}
$\g(C(A))$ is a ring in $K$ such that
$\Z[\la]\subset\g(C(A))\subset\mathcal O_K$, and
$\g(Z(A))=\mathcal U_K\cap\g(C(A))$.
\end{lemma}
\begin{proof}$\g(C(A))$ is a ring because $C(A)$ is a ring.
As we pointed out above images of integer matrices are
algebraic integers and images of matrices with determinant
$\pm 1$ are algebraic units. Hence
$\g(C(A))\subset\mathcal O_K$. Finally, for every polynomial
$p$ with integer coefficients, $p(A)$ is an integer matrix,
hence
$\Z[\la]\subset\g(C(A))$.
\end{proof}

Notice that $\Z(\la)$ is a finite index subring of
$\mathcal O_K$; hence $\g(C(A))$ has the same property.
\begin{remark}
The groups of units in two different rings, say $\mathcal 
O_1\subset\mathcal O_2$,
may coincide. Examples can be found in the table of totally real 
cubic fields in  \cite{C}.
\end{remark}

\begin{prop} \label{Dirichlet} $Z(A)$ is isomorphic to
$\Z^{r_1+r_2-1}\times F$ where $r_1$ is the number the real
embeddings, $r_2$ is the number of pairs of complex conjugate
embeddings of the field
$K$ into $\C$, and $F$ is a finite cyclic group.
\end{prop}
\begin{proof}By lemma \ref{l:rings-units}, $Z(A)$ is isomorphic
to the group of units in the order $\mathcal O$, the statement
follows from the Dirichlet Unit Theorem (\cite{BS}, Ch.2, \S
3).\end{proof}

Now consider an irreducible action $\al$ of $\Z^d$ on $\T^n$.
Denote $\rho_{\al}(\Z^d)$ by $\Gamma$, and
let $\la$ be an eigenvalue of a matrix $A\in\Gamma$ with an
irreducible characteristic polynomial. The centralizers of
$\Gamma$ in
$M(n,\Z)$ and
$GL(n,\Z)$ coincide with $C(A)$ and $Z(A)$ correspondingly. The
field
$K=\Q(\la)$ has degree
$n$ and we can consider the map $\g$ as above. By Lemma
\ref{l:rings-units} $\g(\Gamma)\subset \mathcal U_K$.

For the purposes of purely algebraic considerations in this and
the next section it is convenient to consider actions of
integer $n\times n$ matrices on $\Q^n$ rather than on $\R^n$
and correspondingly to think of $\al$ as an action by automorphisms
of the rational torus $\T^n_{\Q}=\Q^n/\Z^n$.

Let $v=(v_1,\dots,v_n)$ be an eigenvector of $A$ with eigenvalue
$\la$ whose coordinates belong to $K$. Consider the
``projection'' $\pi:\Q^n\to K$ defined by
$\pi(r_1,\dots\,r_n)=\sum_{i=1}^n r_iv_i$. It is a bijection
(\cite{W}, Prop. 8) which conjugates the action of the group
$\Gamma$ with the action on $K$ given by multiplication by
corresponding eigenvalues
$\prod_{i=1}^d\la_i^{k_i},\,\,k_1,\dots,k_d\in\Z$. Here
$A_1,\dots,A_d\in
\Gamma$ are the images of the generators of the action $\al$, and
$A_iv=\la_i v,\,\,i=1,\dots, d$. The lattice $\pi\Z^n\subset K$
is a module over the ring $\Z[\la_1,\dots,\la_d]$.

Conversely, any such data, consisting of an algebraic number
field $K=\Q(\la)$ of degree $n$, a $d$-tuple
$\bar\la=(\la_1,\dots,\la_d)$ of multiplicatively independent
units in $K$, and a lattice $\mathcal L\subset K$
which is a module over $\Z[\la_1,\dots,\la_d]$, determine an
  $\Z^d$-action $\al_{\bar\la,\mathcal
L}$ by automorphisms of $\mathbb{T}^n$ up to
algebraic conjugacy (corresponding to a choice of a basis
in the lattice $\mathcal L$).
This action is generated by multiplications by
$\la_1,\dots,\la_d$ (which preserve $\mathcal L$ by assumption).
The action $\al_{\bar\la,\mathcal L}$ diagonalizes over $\C$ as
follows. Let
$\phi_1=\id,\phi_2,\dots,\phi_n$ be different embeddings of $K$
into $\C$. The multiplications by $\la_i,\,\,i=1,\dots,d$, are
simultaneously conjugate over $\C$ to the respective matrices
\[
\left(\begin{smallmatrix}\la_i & 0 & \dots & 0 \\
0 & \phi_2(\la_i) & \dots & 0 \\
\dots & \dots & \dots & \dots\\
0 & 0 & \dots & \phi_n(\la_i)
\end{smallmatrix}\right),\quad i=1,\dots,d.
\]

We will assume that the action is irreducible  which in many interesting cases
can be easily checked.

Thus, all actions $\al_{\bar\la,\mathcal L}$ with fixed
$\bar\la$ are weakly algebraically isomorphic
since the corresponding embeddings are conjugate over
$\Q$ (Proposition \ref{prop-finfactors}). Actions produced with
different sets of units in the same field, say $\bar{\la}$ and
$\bar{\mu}=(\mu_1,\dots,\mu_d)$, are weakly
algebraically isomorphic if and only if there is an
element $g$ of the Galois group of $K$ such that
$\mu_i=g\la_i,\,\, i=1,\dots,d$. By Proposition
\ref{prop:cyclic} there is a unique cyclic action (up
to algebraic isomorphism) within any class of weakly
algebraically isomorphic actions: it corresponds to setting
$\mathcal L=\Z[\la_1,\dots,\la_d]$; we will
denote this action by
$\al_{\bar\la}^{\min}$. Cyclicity of the action
$\al_{\bar\la}^{\min}$ is obvious since the whole lattice is
obtained from its single element $1$ by the action of the ring
$\Z[\lambda _1^{\pm1},\dots,\lambda _d^{\pm1}]$.

Let us summarize this discussion.

\begin{prop}\label{prop-class} Any irreducible action $\al$ of
$\Z^d$ by automorphisms of $\T^n$ is algebraically conjugate to
an action of the form
   $\al_{\bar\la,\mathcal L}$. It is weakly algebraically
conjugate to the cyclic action $\al_{\bar\la}^{\min}$. The
field $K=\Q[\la_1,\dots,\la_d]$ has degree $n$, and the vector
of units
$\bar\la=(\la_1,\dots,\la_d)$ is defined up to action by an
element of the Galois group of $K:\mathbb{Q}$.
\end{prop}

Apart from the cyclic model $\al_{\bar\la}^{\min}$ there is
another canonical choice of the lattice $\mathcal L$, namely
the ring of integers $\mathcal O_K$. We will denote the action
$\al_{\bar\la,\mathcal O_K}$ by $\al_{\bar\la}^{\max}$.
More generally, one can choose as the lattice $\mathcal L$
any subring
$\mathcal O$ such that $\Z[\la_1,\dots,\la_d]\subset\mathcal
O\subset\mathcal O_K$.

\begin{prop}\label{prop:min-max} Assume that $\mathcal O
\supsetneq \Z[\la_1,\dots,\la_d]$. Then the action
$\al_{\bar{\la},\mathcal O}$ is not algebraically
isomorphic up to a time change to $\al_{\bar\la}^{\min}$. In
particular, if
$\mathcal O_K\neq\Z[\la_1,\dots,\la_d]$, then the actions
$\al_{\bar\la}^{\max}$ and
$\al_{\bar\la}^{\min}$ are not algebraically isomorphic up to
a time change.
\end{prop}

\begin{proof} Let us denote the centralizers in $M(n,\Z)$
of the actions $\al_{\bar{\la},\mathcal
O}$ and $\al_{\bar\la}^{\min}$ by $C_1$ and $C_2$, respectively.
The centralizer $C_1$ contains multiplications by
all elements of $\mathcal O$. For, if one takes any basis
in $\mathcal O$, the multiplication by an element
$\mu\in\mathcal O$ takes elements of the basis into elements of
$\mathcal O$, which are linear combinations with integral
coefficients of the basis elements; hence the multiplication is
given by an integer matrix. On the other hand any element of
each centralizer is a multiplication by an integer in $K$
(Lemma \ref{l:rings-units}).

Now assume that the multiplication by
$\mu\in\mathcal O_K$ belongs to
$C_2$. This means that this multiplication preserves
$\Z[\la_1,\dots,\la_d]$; in particular, $\mu=\mu\cdot 1\in
\Z[\la_1,\dots,\la_d]$. Thus $C_2$ consists of
multiplication by elements of
$\Z[\la_1,\dots,\la_d]$. An algebraic isomorphism up to a
time change has to preserve both the module of polynomials
with integer coefficients in
the generators of the action and the
centralizer of the action in $M(n,\Z)$, which is impossible.
\end{proof}

The central question which appears in connection with our
examples is the classification of weakly algebraically
isomorphic Cartan actions up to algebraic isomorphism.

Proposition \ref{prop:min-max} is useful in distinguishing
weakly algebraically isomorphic actions when $\mathcal
O_K\neq\Z[\la_1,\dots,\la_d]$. Cyclicity also can serve as
a distinguishing invariant.

\begin{corollary}\label{cor-cyclic}The action
$\al_{\bar{\la},\mathcal O}$ is cyclic if and only if
$\mathcal O =\Z[\la_1,\dots,\linebreak[0]\la_d]$.
\end{corollary}

\begin{proof}The action $\al_{\bar\la}^{\min}$ corresponding to
the ring $\Z[\la_1,\dots,\la_d]$ is cyclic by definition
since the ring coincides with the orbit of $1$. By
Proposition \ref{prop:cyclic}, if $\al_{\bar{\la},\mathcal
O}$ were cyclic, it would be algebraically
conjugate to $\al_{\bar\la}^{\min}$, which, by Proposition
\ref{prop:min-max}, implies that $\mathcal O =\Z[\la_1,\dots,\la_d]$.
\end{proof}

The property common to all actions of the
$\al_{\bar{\la},\mathcal O}$ is transitivity of the action of
the centralizer $C(\al_{\bar{\la},\mathcal O})$ on the lattice.
Similarly to cyclicity this property is obviously an invariant
of algebraic conjugacy up to a time change.

\begin{prop}\label{prop-trans}
Any irreducible action $\al$ of $\Z^d$ by automorphisms of
$\T^n$ whose centralizer $C(\al)$ in $M(n,\Z)$ acts
transitively on
$\Z^n$ is algebraically isomorphic to an action
$\al_{\bar{\la},\mathcal O}$, where $\mathcal O\subset\mathcal
O_K$ is a ring which contains $\Z[\la_1,\dots,\la_{d}]$.
\end{prop}

\begin{proof} By Proposition \ref{prop-class} any irreducible
action $\al$ of
$\Z^d$ by automorphisms of $\T^n$ is algebraically conjugate to
an action of the form
$\al_{\bar\la,\mathcal L}$ for a lattice $\mathcal L\subset K$.
Let $C$ be the centralizer of $\al_{\bar\la,\mathcal L}$
in the semigroup of linear endomorphisms of
$\mathcal{L}$.
We fix an element $\beta\in\mathcal L$ with
$C(\alpha )\beta =\mathcal{L}$ and consider conjugation
of the action $\al_{\bar\la,\mathcal L}$ by multiplication by
$\beta^{-1}$; this is simply $\al_{\bar\la,\beta^{-1}\mathcal
L}$. The centralizer of $\al_{\bar\la,\beta^{-1}\mathcal
L}$ acts on the element $1\in\beta^{-1}\mathcal L$ transitively.
By Lemma \ref{l:rings-units} the centralizer consists of all
multiplications by elements of a certain subring
$\mathcal O\subset \mathcal O_K$ which contains
$\Z[\la_1,\dots,\la_{d}]$. Thus $1\in\beta^{-1}\mathcal
L=\mathcal O$.
\end{proof}

\subsection{Structure of algebraic and affine centralizers for
irreducible actions} By Lemma \ref{l:rings-units}, the
centralizer $C(\al)$, as an additive group, is isomorphic to
$\Z^n$ and has an additional ring structure. In the terminology of
Proposition \ref{Dirichlet},
the centralizer $Z(\al)$ for an irreducible
action $\al$ by toral automorphisms is isomorphic to
$\Z^{r_1+r_2-1}\times F$.

An irreducible action $\al$ has {\em maximal rank}
if $d=r_1+r_2-1$. In this case $Z(\al)$ is a finite
extension of $\al$.

Notice that any affine map commuting with
an action $\al$ by toral automorphisms preserves the set
$\Fix(\al)$ of fixed points of the action. This set is always a
subgroup of the torus and hence, for an irreducible action,
always finite. The translation by any element of $\Fix(\al)$
commutes with $\al$ and thus belongs to $Z_\mathit{Aff}(\al)$.
Furthermore, the affine centralizers
$Z_\mathit{Aff}(\al)$ and $C_\mathit{Aff}(\al)$ are generated by these
translations and, respectively,
$Z(\al)$ and
$C(\al)$.

\begin{remark}
Most of the material of this section extends to general
irreducible actions of $\Z^d$ by automorphisms of compact
connected abelian groups; a group possessing such an action
must be a torus or a solenoid (\cite{Sch2, Sch3}). In
the solenoid case, which includes natural extensions of $\Z^d$--actions
by toral endomorphisms, the algebraic numbers
$\la_1,\dots ,\la_d$ which appear in the constructions are not
in general integers. As we mentioned in the
introduction we restrict our algebraic setting
here since we are able to exhibit some of the most interesting and
striking new phenomena using Cartan actions and certain actions
directly derived from them. However, other interesting
examples appear for actions on the torus connected
with not totally real algebraic number fields, actions on solenoids,
and actions on zero-dimensional abelian groups (cf. e.g.
\cite{KiS, Sch1, Sch2, Sch3}).

One can also extend the setup of this section to certain
classes of reducible actions. Since some of these satisfy
condition $(\mathcal R)$ basic rigidity results still hold and
a number of further interesting examples can be constructed.
\end{remark}

\section {Cartan actions} \label{s:ca}

\subsection{Structure of Cartan actions}

Of particular interest
for our study are abelian groups of ergodic automorphisms of
$\T^n$ of maximal possible rank $n-1$ (in agreement with the real
rank of the Lie group $\SnR$).

\begin{definition}An action of $\Z^{n-1}$ on $\T^n$ for $n\ge
3$ by ergodic automorphisms is called a {\em Cartan action}.
\end{definition}

\begin{prop}\label{cartan-str}
Let $\al$ be a Cartan action on $\T^n$.
\begin{enumerate}
\item Any element of $\rho_\al$ other than
identity has real eigenvalues and is hyperbolic and thus Bernoulli.
\item $\al$ is irreducible.
\item The centralizer of $Z(\al)$ is a finite
extension of $\rho_\al(\Z^{n-1})$.
\end{enumerate}
\end{prop}

\begin{proof}
First, let us point out that it is sufficient to prove the
proposition for irreducible actions. For, if $\al$ is not
irreducible, it has a nontrivial irreducible algebraic factor
of dimension, say, $m\le n-1$. Since every factor of an
ergodic automorphism is
ergodic, we thus obtain an action of $\Z^{n-1}$ in $\T^m$ by
ergodic automorphisms. By considering a restriction of this action
to a subgroup of rank $m-1$ which contains an irreducible
matrix, we obtain a Cartan action on $\T^m$. By Statement
3. for irreducible actions, the centralizer of this Cartan
action is a finite extension of $\Z^{m-1}$, and thus cannot
contain $\Z^{n-1}$, a contradiction.

Now assuming that $\al$ is irreducible, take a
   matrix $A\in \rho_\al(\Z^{n-1})$ with
irreducible characteristic polynomial $f$. Such a matrix
exists by Proposition \ref{Berend}. It has
distinct eigenvalues, say
$\la=\la_1,\dots,\la_n$. Consider the correspondence $\g$
defined in \eqref{eq:gamma}.
By Lemma \ref{l:rings-units} for every
$B\in\rho_\al(\Z^{n-1})$ we have $\g(B)\in\mathcal U_K$, hence
the group of units $\mathcal U_K$ in $K$ contains a subgroup
isomorphic to
$\Z^{n-1}$. By the Dirichlet Unit Theorem the rank of the
group of units in $K$ is equal to ${r_1+r_2-1}$, where $r_1$ is
the number of real embeddings and $r_2$ is the number of pairs
of complex conjugate embeddings of
$K$ into $\C$. Since $r_1+2r_2=n$ we deduce that $r_2=0$, so
the field $K$ is totally real, that is all eigenvalues of $A$,
and hence of any matrix in $\rho_\al(\Z^{n-1})$, are real. The
same argument gives Statement 3, since any element of the
centralizer of $\rho_\al(\Z^{n-1})$ in $\GnZ$ corresponds to a
unit in $K$. Hyperbolicity of matrices in $\rho_\al(\Z^{n-1})$
is proved in the same way as Lemma \ref{l:inj}.
\end{proof}

\begin{lemma}\label{hyp} Let $A$ be a hyperbolic matrix  in
$\SnZ$ with irreducible characteristic polynomial and distinct
real eigenvalues. Then every element of the centralizer $Z(A)$
other than $\{\pm 1\}$ is hyperbolic.
\end{lemma}

\begin{proof}Assume that $B\in Z(A)$ is not hyperbolic. As
$B$ is simultaneously diagonalizable with $A$ and has real
eigenvalues, it has an eigenvalue
$+1$ or $-1$. The corresponding eigenspace is rational and
$A$--invariant. Since $A$ is irreducible, this eigenspace has
to coincide with the whole space and hence $B=\pm 1$.
\end{proof}

\begin{corollary} Cartan actions are exactly the maximal rank
irreducible actions corresponding to totally real number fields.
\end{corollary}

\begin{corollary} The centralizer $Z(\al)$ for a Cartan action
$\al$ is isomorphic to $\Z^{n-1}\times\{\pm 1\}$.
\end{corollary}

We will call a Cartan action $\al$ {\em maximal}\, if $\al$ is
an index two subgroup in $Z(\al)$.

Let us point out that $Z_\mathit{Aff}(\al)$ is isomorphic
$Z(\al)\times\Fix(\al)$. Thus, the factor of $Z_\mathit{Aff}(\al)$ by
the subgroup of finite order elements is always isomorphic to
$\Z^{n-1}$. If $\al$ is maximal, this factor is identified with
$\al$ itself. In the next Section we will show (Corollary
\ref{cor-centr}) that for a Cartan action $\al$ on $\T^n,\,\,n\ge 3$ 
the isomorphism type of the
pair
$(Z_\mathit{Aff}(\al),\al)$ is an invariant of the measurable
isomorphism. Thus, in particular, for a maximal Cartan action
the order of the group
$\Fix(\al)$ is a measurable invariant.

\begin{remark} An important geometric distinction between
Cartan actions and general irreducible actions by hyperbolic
automorphisms is the absence of multiple Lyapunov exponents. This
greatly simplifies proofs of various rigidity properties both in
the differentiable and measurable context.
\end{remark}

\subsection{ Algebraically nonisomorphic maximal Cartan
actions}\label{ss:LM}

In Section \ref{ss:irr-integers} we described a particular
class of irreducible actions $\al_{\bar\la,\mathcal O}$
which is characterized by the transitivity of the action of the
centralizer
   $C(\al_{\bar\la,\mathcal O})$ on the lattice
(Proposition \ref{prop-trans}). In the case $\mathcal
O_K=\Z[\la]$ there is only one such action, namely the cyclic
one (Corollary \ref{cor-cyclic}). Now we will analyze this
special case for totally real fields in detail and show how
information about the class number of the field helps to
construct algebraically nonisomorphic maximal Cartan actions.
This will in particular provide examples of Cartan actions not
isomorphic up to a time change to any action of the form
$\al_{\bar\la,\mathcal O}$.

It is well--known that for $n=2$ there are natural bijections
between conjugacy classes of hyperbolic elements in $SL(2,\Z)$
of a given trace, ideal classes in the corresponding real
quadratic field, and congruence classes of primitive integral
indefinite quadratic forms of the corresponding discriminant.
This has been used by Sarnak \cite{Sa} in his proof of the
Prime Geodesic Theorem for surfaces of constant negative
curvature (see also \cite{K}). It follows from an old Theorem
of Latimer and MacDuffee (see \cite{LM},
\cite{T}, and a more modern account in \cite{W}), that the
first bijection persists for $n>2$. Let $A$ a hyperbolic matrix
$A\in \SnZ$ with irreducible characteristic polynomial $f$,
and hence distinct real eigenvalues,
$K=\Q(\la)$, where $\la$ is an eigenvalue of $A$, and $\mathcal
O_K=\Z[\la]$. To each matrix $A'$ with the same eigenvalues,
we assign the eigenvector $v=(v_1,\dots,v_n)$ with
eigenvalue
$\la$: $A'v=\la v$ with all its entries in $\mathcal O_K$,
which can be always done, and to this eigenvector, an ideal in
$\mathcal O_K$ with the
$\Z$--basis
$v_1,\dots,v_n$. The described map is a bijection between the
$GL(n,\Z)$--conjugacy classes of matrices in
$\SnZ$ which have the same characteristic polynomial $f$ and
the set of ideal classes in $\mathcal O_K$. Moreover, it
allows us to reach conclusions about centralizers as well.

\begin{theorem}\label{mytheor} Let $A\in \SnZ$ be a hyperbolic
matrix with irreducible characteristic polynomial $f$
and distinct real eigenvalues,
$K=\Q(\la)$ where $\la$ is an eigenvalue of $A$, and $\mathcal
O_K=\Z[\la]$. Suppose the number of eigenvalues among
$\la_1,\dots,\la_n$ that belong to $K$ is equal to $r$.
If the class number
$h(K)>r$, then there exists a
matrix
$A'\in\SnZ$ having the same eigenvalues as $A$ whose
centralizer $Z(A')$ is not conjugate in $GL(n,\Z)$ to $Z(A)$.
Furthermore, the number of matrices in $\SnZ$ having the same
eigenvalues as $A$ with pairwise nonconjugate \textup{(}in
$GL(n,\Z)$\textup{)} centralizers is at least $[\frac{h(K)}{r}] +1$, where
$[x]$ is the largest integer $<x$.
\end{theorem}

\begin{proof}Suppose the matrix $A$ corresponds to the ideal
class $I_1$ with the $\Z$--basis $v^{(1)}$. Then
\[ A v^{(1)}=\la v^{(1)}.\] Since $h(K)>1$, there exists a
matrix $A_2$ having the same eigenvalues which corresponds to a
different ideal class $I_2$ with the basis
$v^{(2)}$, and we have
   \[ A_2 v^{(2)}=\la v^{(2)}. \] The eigenvectors $v^{(1)}$ and
$v^{(2)}$ are chosen with all their entries in $\mathcal O_K$.
Now assume that $Z(A_2)$ is conjugate to $Z(A)$. Then
$Z(A_2)$ contains a matrix $B_2$ conjugate to $A$. Since $B_2$
commutes with
$A_2$ we have
$B_2 v^{(2)}=\mu_2 v^{(2)}$, and since $B_2$ is conjugate to
$A$, $\mu_2$ is one of the roots of $f$. Moreover, since
$B_2\in \SnZ$ and all entries of
$v^{(2)}$ are in $K$, $\mu_2\in K$. Thus $\mu_2$ is one of $r$
roots of $f$ which belongs to $K$.

  From $B_2=S^{-1} A S$
($S\in GL(n,\Z)$) we deduce that $\mu_2
(Sv^{(2)})=A(Sv^{(2)})$.
Since $I_1$ and $I_2$ belong to different ideal classes,
$Sv^{(2)}\ne kv^{(1)}$ for any $k$ in the quotient field of
$\mathcal O_K$, and since $\la$ is
a simple eigenvalue for $A$, we deduce that
$\mu_2\ne\la$, and thus $\mu_2$ can take one of the $r-1$
remaining values.

Now assume that
$A_3$ corresponds to the third ideal class, i.e
\[
A_3 v^{(3)}=\la v^{(3)},
\]
and $B_3$ commutes
with
$A_3$ and is conjugate to $A$, and hence to $B_2$. Then
$B_3 v^{(3)}=\mu_3 v^{(3)}$ where $\mu_3$ is a root of $f$
belonging to the field $K$.
By the
previous considerations,
$\mu_3\ne\la$ and $\mu_3\ne\mu_2$. An induction argument shows
that if the class number of $K$ is greater than $r$, there
exists a matrix $A'$ such that no matrix in
$Z(A')$ is conjugate to $A$, i.e. $Z(A')$ and $Z(A)$ are not
conjugate in
$GL(n,\Z)$.

Since $A'$ has the same characteristic polynomial as $A$,
continuing the same process, we can find not more than $r$
matrices representing different ideal classes having
centralizers conjugate to $Z(A')$, and the required
estimate follows.
\end{proof}

\section {Measure--theoretic rigidity of conjugacies,
centralizers, and factors}\label{s:rigidity}

\subsection{Conjugacies}\label{ss:conj}
Suppose
$\al$ and $\al'$ are measurable actions of the same group $G$ by
measure--preserving transformations of the spaces $(X,\mu)$ and
$(Y,\nu)$, respectively. If $H:(X,\mu)\to(Y,\nu)$ is a metric
isomorphism (conjugacy) between the actions then the lift of
the measure $\mu$ onto the
$\mbox{graph} \,H\subset X\times Y$ coincides with the lift of
$\nu$ to
$\mbox{graph}\, H^{-1}$. The resulting measure $\eta$ is a very
special case of a {\it joining} of $\al$ and $\al'$: it is
invariant under the diagonal (product) action $\al\times\al'$
and its projections to $X$ and
$Y$ coincide with $\mu$ and $\nu$, respectively. Obviously
the projections establish metric isomorphism of the action
$\al\times\al'$ on
$(X\times Y, \eta)$ with $\al$ on $(X,\mu)$ and $\al'$ on
$(Y,\nu)$ correspondingly.

Similarly, if an automorphism $ H:(X,\mu)\to(X,\mu)$ commutes
with the action
$\al$, the lift of $\mu$ to $\mbox{graph}
\,H\subset X\times X$ is a self-joining of $\al$, i.e. it is
$\al\times\al$--invariant and both of its projections coincide
with $\mu$. Thus an information about invariant measures of the
products of different actions as well as the product of an
action with itself may give an information about isomorphisms
and centralizers.

The use of this joining construction in
order to deduce rigidity of isomorphisms and centralizers from
properties of invariant measures of the product was first
suggested in this context to the authors by J.-P. Thouvenot.

In both cases the ergodic properties of the
joining would be known because of the isomorphism with the
original actions. Very similar considerations apply to the
actions of semi--groups by noninvertible measure--preserving
transformations. We will use the following corollary of the
results of \cite{KS3}.

\begin{theorem}\label{corr-K} Let $\al$ be an action of
${\Z^2}$ by ergodic toral automorphisms
and let $\mu$ be a weakly mixing
$\al$--invariant measure such that for some
$\mathbf{m}\in \Z^2$, $\al^\mathbf{m}$ is a $K$-automorphism.
Then $\mu$ is a
translate of Haar measure on an
$\al$--invariant rational subtorus.
\end{theorem}

\begin{proof} We refer to Corollary 5.2' from (\cite{KS3},
``Corrections...'').  According to this corollary the measure $\mu$ is
an extension of a zero entropy measure for an algebraic factor of
smaller dimension  with Haar conditional measures in the fiber. But
since $\al$ contains a  $K$-automorphism it does not have non--trivial
zero entropy factors. Hence the factor in question is the action on a
single point and $\mu$ itself is a Haar measure on a rational
subtorus.  \end{proof}

Conclusion of Theorem~\ref{corr-K}  obviously holds for any action of
$\Zd,\,\,d\ge 2$ which contains a subgroup $\Z^2$ satisfying
assumptions of Theorem~\ref{corr-K}. Thus  we can deduce the following
result which is central for our constructions.

\begin{theorem} \label{thm-isom} Let $\al$ and $\al'$ be two
actions of
$\Zd$ by automorphisms of $\T^n$ and $\T^{n'}$ correspondingly
and assume that $\al$ satisfies condition $(\mathcal R)$.
Suppose that $H:\T^n\to \T^{n'}$ is a measure--preserving
isomorphism between
$(\al, \la)$ and $(\al', \la)$, where $\la$ is Haar measure.
Then
$n=n'$ and $H$ coincides \textup{(mod\;0)} with an affine automorphism
on the torus
$\T^n$, and hence $\al$ and $\al'$ are algebraically
isomorphic. \end{theorem}
\begin{proof} First of all, condition
$(\mathcal R)$ is invariant under metric isomorphism, hence
$\al'$ also satisfies this condition. But ergodicity with
respect to Haar measure can also be expressed in terms of the
eigenvalues; hence $\al\times\al'$ also satisfies ($\mathcal
R$).
   Now consider the joining measure $\eta$ on $\mbox{graph}
\,H\subset
\T^{n+n'}$. The conditions of Theorem~\ref{corr-K} are satisfied
for the invariant measure $\eta$ of the action
$\al\times\al'$. Thus $\eta$ is a translate of Haar measure
on a rational $\al\times\al'$--invariant subtorus
$\T'\subset \T^{n+n'}=\T^n\times \T^{n'}$. On the other hand we
know that projections of $\T'$ to both $\T^n$ and $\T^{n'}$
preserve Haar measure and are one--to--one. The partitions of
$\T'$ into pre--images of points for each of the projections are
measurable partitions and Haar measures on elements are
conditional measures. This implies that both projections are
onto, both partitions are partitions into points, and hence
$n=n'$ and
$\T'=\mbox{graph}\, I$, where $I: \T^n\to \T^n$ is an affine
automorphism which has to coincide $(\textup{mod}\;0)$ with the
measure--preserving isomorphism
$H$. \end{proof}

Since a time change is in a sense a trivial modification of an
action we are primarily interested in distinguishing actions up
to a time change. The corresponding rigidity criterion follows
immediately from Theorem \ref{thm-isom}.

\begin{corollary}\label{cor:time-change} Let $\al$ and $\al'$ be
two actions of
$\Zd$ by automorphisms of $\T^n$ and $\T^{n'}$, respectively,
and assume that
$\al$ satisfies condition $(\mathcal R)$. If $\al$ and $\al'$
are measurably isomorphic up to a time change then they are
algebraically isomorphic up to a time change. \end{corollary}

\subsection {Centralizers} Applying Theorem~\ref{thm-isom} to
the case
$\al=\al'$ we immediately obtain rigidity of the centralizers.

\begin{corollary} \label{cor-centr} Let
$\al$ be an action of $\Zd$ by automorphisms of $\T^n$
satisfying condition $(\mathcal R)$. Any invertible Lebesgue
measure--preserving transformation commuting with $\al$
coincides \textup{(mod 0)} with an affine automorphism of
$\T^n$. \end{corollary}

Any affine transformation commuting with $\al$ preserves the
finite set
of fixed points of the action. Hence the centralizer of $\al$
in affine automorphisms has a finite index subgroups
which consist of automorphisms and which corresponds to the
centralizer of
$\rho_{\al}(\Z^d)$ in $GL(n,\Z)$.

Thus, in contrast with the case of a single automorphism, the
centralizer of such an action $\al$ is not more than countable,
and can be identified with a finite extension of a certain
subgroup of $GL(n,\Z)$. As an immediate consequence we obtain
the following result.

\begin{prop} For any $d$ and $k$, $2\le d\le k$, there
exists a
$\Zd$--action by hyperbolic toral automorphisms such that its
centralizer in the group of Lebesgue measure--preserving
transformations is isomorphic to
$\{\pm 1\}\times\Z^k$.
\end{prop}

\begin{proof} Consider a hyperbolic matrix $A\in
SL(k+1,\Z)$ with irreducible characteristic polynomial and
real eigenvalues such that the origin is the only fixed point
of $F_A$. Consider a subgroup of $Z(A)$
isomorphic to $\Zd$ and containing $A$ as one of its
generators. This subgroup determines an embedding $\rho:\Z^d\to
SL(k+1,\Z)$. Since $d\ge 2$ and by Proposition \ref{hyp},
all matrices in
$\rho(\Zd)$ are hyperbolic and hence ergodic, condition
$(\mathcal R)$ is satisfied. Hence by Corollary
\ref{cor-centr}, the measure--theoretic centralizer of the
action $\al_{\rho}$ coincides with its algebraic centralizer,
which, in turn, and obviously, coincides with centralizer of the
single automorphism $F_A$ isomorphic to
$\{\pm 1\}\times\Z^k$.
\end{proof}

\subsection{Factors, noninvertible centralizers
and weak isomorphism}

A small modification of the proof of Theorem \ref{thm-isom}
produces a result about rigidity of factors.

\begin{theorem} \label{thm-factors} Let $\al$ and $\al'$ be two
actions of
$\Zd$ by automorphisms of $\T^n$ and $\T^{n'}$
respectively, and assume that $\al$ satisfies condition
$(\mathcal R)$. Suppose that $H:\T^n\to \T^{n'}$ is a
Lebesgue measure--preserving transformation such that
$H\circ\al=\al'\circ H$.
   Then $\al'$ also satisfies $(\mathcal R)$
   and $H$ coincides \textup{(mod 0)} with an epimorphism $h:\T^n\to
\T^{n'}$ followed by translation. In particular, $\al'$ is an
algebraic factor of $\al$.
   \end{theorem}

\begin{proof} Since $\al'$ is a measurable factor of
$\al$, every element which is ergodic for $\al$ is also ergodic
for $\al'$. Hence $\al'$ also satisfies condition $(\mathcal
R)$.
   As before consider the product action $\al\times\al'$ which
now by the same argument also satisfies $(\mathcal R)$. Take
the $\al\times\al'$ invariant measure $\eta=(\mbox{Id}\times H)
_*\la$ on $\mbox{graph} \,H$. This measure provides a joining
of $\al$ and $\al'$. Since $(\al\times\al',
(\mbox{Id}\times H)_*\la)$ is isomorphic to $(\al,\la)$ the
conditions of Corollary \ref{corr-K} are satisfied and $\eta$
is a translate of Haar measure on an invariant rational
subtorus
$\T'$. Since $\T'$ projects to the first coordinate one-to-one we
deduce that
$H$ is an algebraic epimorphism \textup{(mod 0)} followed by a
translation.
\end{proof}

Similarly to the previous section the application of
Theorem \ref{thm-factors}
to the case
$\al=\al'$ gives a description of the centralizer of $\al$
in the group of all measure--preserving transformations.

\begin{corollary} \label{cor-centrNI} Let
$\al$ be an action of $\Zd$ by automorphisms of $\T^n$
satisfying condition $(\mathcal R)$. Any Lebesgue
measure--preserving transformation commuting with $\al$
coincides \textup{(mod 0)} with an affine map on
$\T^n$. \end{corollary}

Now we can obtain the following strengthening of Proposition
\ref{prop-finfactors}
   for actions satisfying condition $(\mathcal R)$
which is one of the central conclusions of this paper.

\begin{theorem}\label{thm:weak-iso} Let $\al$ be an action of
$\Z^d$ by automorphisms of $\T^n$ satisfying condition
$(\mathcal R)$ and
$\al'$ another $\Z^d$-action by toral automorphisms. Then
$(\al,\,\la)$ is weakly isomorphic to $(\al',\,\la')$ if and
only if $\rho_{\al}$ and $\rho_{\al'}$ are isomorphic over
$\Q$, i.e. if
$\al$ and $\al'$ are finite algebraic factors of each other.
\end{theorem}

\begin{proof} By Theorem~\ref{thm-factors}, $\al$ and $\al'$ are
algebraic factors of each other. This implies that $\al'$ acts
on the torus of the same dimension $n$ and hence both algebraic
factor--maps have finite fibres. Now the statement follows
from Proposition \ref{prop-finfactors}.
\end{proof}

\subsection{Distinguishing weakly isomorphic actions}

Similarly we can translate criteria for algebraic
conjugacy of weakly algebraically conjugate actions to the
measurable setting.

\begin{theorem}\label{thm:cyclic}If $\al$ is an irreducible
cyclic action of $\Z^d,\,\,d\ge 2$, on $\T^n$ and $\al'$ is a
non--cyclic
$\Z^d$-action by toral automorphisms. Then $\al$ and $\al'$
are not measurably isomorphic up to a time change.
\end{theorem}

\begin{proof} Since action $\al$ satisfies condition
$(\mathcal R)$ (Corollary \ref{cor:irred}) we can apply
Theorem \ref {thm:weak-iso} and conclude that we only need to
consider the case when $\rho_{\al}$ and $\rho_{\al'}$ are
isomorphic over $\Q$ up to a time change. But then, by
Proposition \ref{prop:cyclic}, $\al$ and $\al'$ are not
algebraically isomorphic up to a time change and hence, by
Corollary
\ref{cor:time-change}, they are not measurably isomorphic
up to a time change.
\end{proof}

Combining Proposition \ref{prop:min-max} and Corollary
\ref{cor:time-change} we immediately obtain rigidity for the
minimal irreducible models.

\begin{corollary}\label{cor:min-max} Assume that $\mathcal O\supsetneq
\Z[\la_1,\dots,\la_d]$. Then the action
$\al_{\bar\la,\mathcal O}$ is not measurably
isomorphic up to a time change to
$\al_{\bar\la}^{\min}$. In particular, if $\mathcal
O_K\supsetneq\Z[\la_1,\dots,\la_d]$, then the actions
$\al_{\bar\la}^{\max}$ and $\al_{\bar\la}^{\min}$ are not
measurably isomorphic up to a time change. \end{corollary}

\section {Examples}\label{s:examples}

Now we proceed to produce examples of actions for which the
entropy data coincide but which are not algebraically
isomorphic, and hence by Theorem~\ref{thm-isom} not
measure--theoretically isomorphic.

\subsection{Weakly nonisomorphic actions}\label{ss:not
weakly} In this section we consider
actions which are not algebraically isomorphic over $\Q$ (or,
equivalently, over $\R$) and hence by Theorem \ref
{thm:weak-iso} are not even weakly isomorphic. The easiest way
is as follows.

\medskip

\noindent{\bf Example 1a.} Start with any action
$\al$ of
$\Zd,\,d\ge 2$, by ergodic automorphisms of $\T^n$. We may
double the entropies of all its elements in two different ways:
by considering the Cartesian square $\al\times\al$ acting on
$\T^{2n}$, and by taking second powers of all elements:
$\al_2^\bn=\al^{2\bn}$ for all $\bn\in\Z^d$. Obviously
$\al\times\al$ is not algebraically isomorphic to $\al_2$,
since, for example, they act on tori of different dimension.
Hence by Theorem~\ref{thm-isom}
$(\al\times\al, \la)$ is not metrically isomorphic to
$(\al_2, \la)$ either.

Now we assume that $\al$ contains an automorphism $F_A$ where
$A$ is hyperbolic with an irreducible characteristic polynomial
and distinct positive real eigenvalues. In
this case it is easy to find an invariant distinguishing the two
actions, namely, the algebraic type of the centralizer of the
action in the group of measure--preserving transformations. By
Corollary~\ref{cor-centr}, the centralizer of $\al$ in the
group of measure--preserving transformations coincides with the
centralizer in the group of affine maps,
which is a finite extension of the
centralizer in the group of automorphisms.
By the Dirichlet Unit Theorem, the centralizer of
$Z(\al_2)$ in the group of automorphisms of the torus is isomorphic to
$\{\pm 1\}\times\Z^{n-1}$, whereas the centralizer of
$\al \times\al$ contains the $\Z^{2(n-1)}$--action by product
transformations
$\al^{\mathbf{n}_1}\times\al^{\mathbf{n}_2},\,\,\mathbf{n}_1,
\mathbf{n}_2 \in\Z^{n-1}$.
In fact, the centralizer of $\al\times\al$ can be calculated explicitly:

\begin{prop}Let $\la$ be an eigenvalue of $A$. Then $K=\Q(\la)$
is a totally real algebraic field. If its ring of
integers $\mathcal O_K$ is equal to $\Z[\la]$ then the centralizer of
$\al\times\al$ in $GL(2n,\Z)$ is isomorphic to the group
$GL(2,\mathcal O_K)$, i.e. the group of $2\times 2$ matrices
with entries in
$\mathcal O_K$ whose determinant
is a unit in $\mathcal O_K$.
\end{prop}

\begin{proof}First we notice that a matrix in block form
$B=\left(
	\begin{smallmatrix} X & Y\\Z & T\end{smallmatrix}\right)$
with $X,Y,Z,T\in M(n,\Z)$ commutes with $\left(
	\begin{smallmatrix} A & 0\\0 & A\end{smallmatrix}\right)$ if
an only if
$X,Y,Z,T$ commute with $A$ and can thus be identified with
elements of $\mathcal O_K$. In this case $B$ can be identified with a
matrix in
$M(2,\mathcal O_K)$. Since $\det\left(
	\begin{smallmatrix} X & Y\\Z &
T\end{smallmatrix}\right)=\det(XT-YZ)=\pm 1$ (cf.\ \cite{Ga}), the
norm of the determinant of the $2\times 2$ matrix corresponding
to $B$ is equal $\pm 1$. Hence this determinant is a unit in
$\mathcal O_K$, and we obtain the desired isomorphism.
\end{proof}

It is not difficult to modify Example 1a to obtain weakly
nonisomorphic actions with the same entropy on the torus
of the same dimension.

\medskip

\noindent{\bf Example 1b.}
For a natural number $k$ define the action
$\al_k$ similarly to $\al_2$: $\al_k^{\bn}=\al^{k\bn}$
for all
$\bn\in\Z^d$.

   The actions
$\al_3\times
\al$ and
$\al_2\times\al_2$ act on $\T^{2n}$, have the same
entropies for all elements and are not isomorphic.

As before, we can see that centralizers of these two actions
are not isomorphic. In particular, the centralizer of
$\al_3\times \al$ is abelian since it has simple
eigenvalues, while the centralizer of
$\al_2\times\al_2$ is not.

\subsection{Cartan actions distinguished by
  cyclicity or maximality}\label{ss:cent} We give two examples
which illustrate the method of Section \ref{ss:irr-integers}.
They provide weakly algebraically isomorphic Cartan actions of
$\Z^2$ on
$\T^3$ which are not algebraically isomorphic even up to a time
change (i.e. a linear change of coordinates in
$\Z^2$) by
Proposition \ref {prop:min-max}. These examples utilize the
existence of number fields $K=\Q(\la)$  and units  $\bar\la=(\la_1,\la_2)$
in them for which
$\mathcal O_K\ne\Z[\la_1,\la_2]$.
In each example
one action has a form
$\al_{\bar\la}^{\min}$ and the other $\al_{\bar\la}^{\max}$.
Hence by Corollary \ref{cor:min-max} they are not measurably
isomorphic up to a time change

In other words, in each example one action, namely,
$\al_{\bar\la}^{\min}$, is a cyclic Cartan action, and the
other is not.

We will  aslo show that  in
these examples the conjugacy type of the pair
$(Z(\al),\al)$ distinguishes weakly isomorphic actions. Let us
point out that a noncylic action for example
$\al_{\bar\la}^{\max}$ may be maximal, for example when
fundamental units lie in a proper subring of $\mathcal O_K$.
However in our examples centralizers for the cyclic actions
will be  dirrefern and thus will serve as a distuinguishing invariant.

The information
about cubic fields is either taken from \cite{C} or obtained
with the help of the computer package Pari-GP. Some
calculations were made by Arsen Elkin during the REU program
at Penn State in summer of 1999.

We construct
two $\Z^2$--actions, $\al$, generated by commuting
matrices $A$ and $B$, and
$\al'$, generated by commuting matrices $A'$ and $B'$ in
$GL(3,\Z)$. These actions are weakly algebraically isomorphic
by Proposition
\ref{prop-class} since they are produced with the same set of
units on two different orders,
$\Z[\la]$ and $\mathcal O_K$, but not algebraically isomorphic
by Proposition
\ref{prop:min-max}. In these examples the action $\al$ is
cyclic by Corollary \ref{cor-cyclic} and will be shown to be  a maximal
Cartan action. Thus
$Z(\al)=\al\times\{\pm\Id\}$. The action $\al'$ is not maximal,
specifically, $Z(\al')/\{\pm\Id\}$ is a nontrivial finite
extension of $\al'$.

\medskip

\noindent{\bf Example 2a}. Let $K$ be a totally real cubic
field given by the irreducible polynomial
$f(x)=x^3 + 3 x^2 - 6 x + 1$, i.e.
$K=\Q(\la)$ where $\la$ is one of its roots. The discriminant
of $K$ is equal to $81$, hence its Galois group is cyclic, and
$[\mathcal O_K:\Z[\la]]=3$. The algebraic integers $\la_1=\la$
and $\la_2=2-4\la-\la^2$ are units with $f(\la_1)=f(\la_2)=0$.
The minimal order in $K$ containing $\la_1$ and $\la_2$ is
$\Z[\la_1,\la_2]=\Z[\la]$, and the maximal order is $\mathcal
O_K$.
A basis in fundamental units is $\epsilon=\frac{\la^2+5\la+1}{3}$
and $\epsilon-1$, hence $\mathcal U_K$ is not contained in $\Z[\la]$.

With
respect to the basis $\{1,\la,\la^2\}$ in $\Z[\la]$,
multiplications by $\la_1$ and $\la_2$ are given by the matrices
\[
A=\left(
	\begin{smallmatrix}
\hphantom{-}0&\hphantom{-}1&\hphantom{-}0
	\\
\hphantom{-}0&\hphantom{-}0&\hphantom{-}1
	\\
-1&\hphantom{-}6&-3
	\end{smallmatrix}
\right),\qquad
B=\left(
	\begin{smallmatrix}
2&-4&-1
	\\
1&-4&-1
	\\
1&-5&-1
	\end{smallmatrix}
\right),
	\]
respectively (if acting from the right on row--vectors).
A direct calculation shows that this action is maximal.

With respect to the basis $\{-\frac 2{3}+\frac 5{3}\la+\frac
1{3}\la^2, -\frac 1{3}+\frac 7{3}\la+\frac 2{3}\la^2\}$ in
$\mathcal O_K$, multiplications by $\la_1$ and $\la_2$ are given
by the matrices
\[
A'=\left(
	\begin{smallmatrix}
\hphantom{-}1&\hphantom{-}2&-1
	\\
-1&-2&\hphantom{-}2
	\\
\hphantom{-}2&\hphantom{-}5&-2
	\end{smallmatrix}
\right),\qquad
B'=\left(
	\begin{smallmatrix}
\hphantom{-}1&-1&-1
	\\
-1&-2&-1
	\\
-1&-4&-2
	\end{smallmatrix}
\right).
\]
We have $A'=VAV^{-1}$, $B'=VBV^{-1}$ for $
V=\left(
	\begin{smallmatrix}
2&-2&-1
	\\
0&-3&\hphantom{-}0
	\\
1&-4&-2
	\end{smallmatrix}
\right)$.
Since $A$ is a companion matrix of $f$, $\al=\langle
A,B\rangle$ has a cyclic element in $\Z^3$. If $A'$ also had a
cyclic element $\mathbf m=(m_1,m_2,m_3)\in\Z^3$, then the vectors
\[
\begin{smallmatrix}
\mathbf m=(m_1,m_2,m_3),\enspace
\mathbf
mA'=(m_1-m_2+2m_3,2m_1-2m_2+5m_3,-m_1+2m_2-2m_3)\\
\mathbf
m(A')^2=(-3m_1+5m_2-7m_3,-7m_1+12m_2-16m_3,5m_1-7m_2+12m_3),
\end{smallmatrix}
\]
would have to generate $\Z^3$ or, equivalently
	\begin{gather*}
\det \left(
	\begin{smallmatrix}
m_1&m_2&m_3
	\\
m_1 - m_2 + 2 m_3&2 m_1 - 2 m_2 + 5 m_3&-m_1+2m_2-2m_3
	\\
- 3 m_1 + 5 m_2 - 7 m_3&- 7 m_1 + 12 m_2 - 16
m_3&5 m_1 - 7 m_2 +
12 m_3
   	\end{smallmatrix}
\right)
	\\
	\begin{aligned}
= 3 m_1^3 &+ 18 m_1^2m_3 - 9 m_1 m_2^2 - 9 m_1 m_2 m_3
	\\
&+ 27 m_1 m_3^2 + 3 m_2 ^3 - 9 m_2 m_3^2+ 3 m_3^3=1.
	\end{aligned}
	\end{gather*}
This contradiction shows that $A'$ has no cyclic vector, and
since $B'=2-4A'-{A'}^2$ , the action $\al'$ is not cyclic.
In this example both actions $\al$ and $\al'$ have a single
fixed point $(0,0,0)$, hence their linear and affine
centralizers coincide, and by Corollary \ref{cor:time-change}
$\al$ and
$\al'$ are not measurably isomorphic up to a time change.

The action $\al'$ is not maximal beacuse $Z(\al')$
contains fundamental units.
\medskip

\noindent{\bf Example 2b}. Let us consider a totally
real cubic field $K$ given by the irreducible polynomial
$f(x)=x^3 - 7x^2 + 11x - 1$. Thus
$K=\Q(\la)$ where $\la$ is one of its roots. In this field
the ring of integers $\mathcal O_K$ has basis $\{1,\la,\frac
1{2}\la^2+\frac1{2}\}$ and hence $[\mathcal O_K:\Z[\la]]=2$.
The
fundamental units in $\mathcal O_K$ are $\{\frac
1{2}\la^2-2\la+\frac 1{2}, \la-2\}$. We choose the units
$\la = \la_1=(\frac 1{2}\la^2-2\la+\frac 1{2})^2$ and
$\la_2=\la-2$ which are contained in both orders, $\mathcal
O_K$ and $\Z[\la]$.

In $\Z[\la]$ we consider the basis
$\{1,\la,\la^2\}$ relative to which the multiplication by
$\la_1$ is represented by the companion matrix
$A=\left(
	\begin{smallmatrix} \hphantom{-}0 & \hphantom{-}1 &
\hphantom{-}0\\ \hphantom{-}0 & \hphantom{-}0 & \hphantom{-}1\\
\hphantom{-}1 & -11 & \hphantom{-}7\end{smallmatrix}\right)$
and multiplication by
$\la_2$ is represented by the matrix $B=\left(
	\begin{smallmatrix} -2 & \hphantom{-}1 & \hphantom{-}0\\
\hphantom{-}0 & -2 & \hphantom{-}1\\\hphantom{-}1 & -11 &
\hphantom{-}5\end{smallmatrix}\right)$.

For $\mathcal O_K$ with the basis $\{1,\la,\frac
1{2}\la^2+\frac1{2}\}$ multiplications by $\la_1$ and $\la_2$
are represented by the matrices $A'=\left(
	\begin{smallmatrix} \hphantom{-}0 & \hphantom{-}1 &
\hphantom{-}0\\ -1 & \hphantom{-}0 & \hphantom{-}2\\
-3 & -5 & \hphantom{-}7\end{smallmatrix}\right)$ and
$B'=\left(
	\begin{smallmatrix} -2 & \hphantom{-}1 &
\hphantom{-}0\\ -1 & -2 & \hphantom{-}2\\
-3 & -5 & \hphantom{-}5\end{smallmatrix}\right)$.

It can be seen directly that $\al$ and $\al'$ are not
algebraically conjugate up to a time change since $A'$ is a square of a
matrix from
$SL(3,\Z)$:
$A'= {\left(\begin{smallmatrix}\hphantom{-}0 & -2 &
\hphantom{-}1\\ -1 & -5 & \hphantom{-}3\\
-2 & -9 & \hphantom{-}6\end{smallmatrix}\right)}^2$, while $A$ is
not a square of a matrix in
$GL(3,\Z)$, which is checked by reducing modulo $2$. In this
case it is also easily seen that the action $\al'$ is not
cyclic since the corresponding determinant is divisible by 2.
The action $\al$ has $2$ fixed points on $\T^3$: $(0,0,0)$ and
$(\frac 1{2},\frac 1{2},\frac 1{2})$, while the action
$\al'$ has $4$ fixed points: $(0,0,0)$,
$(\frac 1{2},\frac 1{2},\frac 1{2})$, $(\frac 1{2},\frac
1{2},0)$, and $(0,0,\frac 1{2})$. Hence the affine
centralizer of $\al$ is $Z(\al)\times\Z/2\Z$, and the affine
centralizer of $\al'$ is $Z(\al')\times(\Z/2\Z\times\Z/2\Z)$.

   By Lemma
\ref{hyp},
the group of elements of finite order in $Z_\mathit{Aff}(\al)$ is
$\Z/2\Z\times
\Z/2\Z$ and in $Z_\mathit{Aff}(\al')$ it is $\Z/2\Z\times
\Z/2\Z\times \Z/2\Z$. The indices of each
action in its affine centralizer are $[Z_\mathit{Aff}(\al):\al]=4$ and
$[Z_\mathit{Aff}(\al'):\al']=16$.

This gives two alternative arguments
that the actions are not measurably isomorphic up to a time
change.
\medskip

\subsection{Nonisomorphic maximal Cartan
actions}\label{ss:max-Cart} We find examples of weakly
algebraically isomorphic maximal Cartan actions which are not
algebraically isomorphic up to time change. For such an
action $\al$ the structure of the pair $(Z(\al),\al)$ is always
the same: $Z(\al)$ is isomorphic as a group to
$\al\times\{\pm\Id\}$. The algebraic tool which allows to
distinguish the actions is Theorem \ref{mytheor}.

\medskip

\noindent{\bf Example 3a}. An
example for
$n=3$ can be obtained from a totally real cubic field with class
number
$2$ and the Galois group $S_3$. The smallest discriminant for
such a field is $1957$ (\cite{C}, Table B4), and it can be
represented as
$K=\Q(\la)$ where $\la$ is a unit in $K$ with minimal polynomial
$f(x)=x^3-2x^2-8x-1$. In this field the ring of integers
$\mathcal O_K=\Z[\la]$ and the
fundamental units are $\la_1=\la$ and $\la_2=\la+2$. Two actions
are constructed with this set of units (fundamental, hence
multiplicatively independent) on two different lattices,
$\mathcal O_K$ with the basis
$\{1,\la,\la^2\}$, representing the principal ideal class, and
$\mathcal L$ with the basis $\{2, 1+\la, 1+\la^2\}$
representing to the second ideal class. Notice that the units
$\la_1$ and
$\la_2$ do not belong to $\mathcal L$, but $\mathcal L$ is a
$\Z[\la]$-module. The first action $\al$ is generated by the
matrices
$A=\left(
	\begin{smallmatrix} 0 & 1& 0\\ 0 & 0 & 1\\ 1 & 8 & 2
\end{smallmatrix}\right)$ and
$B=\left(	\begin{smallmatrix} 2
& 1 & 0\\ 0 & 2 & 1\\ 1 & 8 & 4
\end{smallmatrix}\right)$ which represent multiplication by
$\la_1$ and $\la_2$, respectively, on $\mathcal O_K$. The
second action $\al'$ is generated by matrices
$A'=\left(
	\begin{smallmatrix} -1 & \hphantom{-}2 & \hphantom{-}0\\ -1 &
\hphantom{-}1 & \hphantom{-}1\\ -5 & \hphantom{-}9 &
\hphantom{-}2\end{smallmatrix}\right)$ and
$B'=\left(
	\begin{smallmatrix} \hphantom{-}1 & \hphantom{-}2 &
\hphantom{-}0\\ -1 &
\hphantom{-}3 & \hphantom{-}1\\ -5 & \hphantom{-}9 &
\hphantom{-}5\end{smallmatrix}\right)$ which represent
multiplication by
$\la_1$ and $\la_2$, respectively, on $\mathcal L$ in the given
basis. By Proposition \ref{prop-class} these actions are weakly
algebraically isomorphic. By Theorem
\ref{mytheor} they are not algebraically isomorphic. Since the
Galois group is $S_3$ there are no nontrivial time changes
which produce conjugacy over $\Q$. Therefore, but Theorem
\ref{thm-isom} the actions are not measurably isomorphic.

It is interesting to point out that for actions $\al$ and $\al'$
the affine centralizers $Z_\mathit{Aff}(\al)$ and $Z_\mathit{Aff}(\al')$ are
not isomorphic as abstract groups. The
action
$\al$ has
$2$ fixed points on
$\T^3$:
$(0,0,0)$ and
$(\frac 1{2},\frac 1{2},\frac 1{2})$, while the action
$\al'$ has a single fixed point $(0,0,0)$. Hence
$Z_\mathit{Aff}(\al)$ is isomorphic to $Z(\al)\times\Z/2\Z$,
$Z_\mathit{Aff}(\al')$ is isomorphic to $Z(\al')$. As abstract
groups, $Z_\mathit{Aff}(\al)\approx \Z^2\times \Z/2\Z\times \Z/2\Z$
and
$Z_\mathit{Aff}(\al')\approx \Z^2\times \Z/2\Z$.

Hence by Corollary
\ref{cor-centr} the measurable centralizers of $\al$ and $\al'$
are not conjugate in the group of measure--preserving
transformation providing a distinguishing invariant of
measurable isomorphism.

\medskip

\noindent{\bf Example 3b}. This example is obtained from a
totally real cubic field with class number $3$,
Galois group $S_3$, and discriminant $2597$. It can be
represented as
$K=\Q(\la)$ where $\la$ is a unit in $K$ with minimal polynomial
$f(x)=x^3-2x^2-8x+1$. In this field the ring of integers
$\mathcal O_K=\Z[\la]$ and the
fundamental units are $\la_1=\la$ and $\la_2=\la+2$. Three
actions are constructed with this set of units on three
different lattices,
$\mathcal O_K$ with the basis
$\{1,\la,\la^2\}$, representing the principal ideal class,
$\mathcal L$ with the basis $\{2, 1+\la, 1+\la^2\}$
representing the second ideal class, and $\mathcal L^2$ with
the basis $\{4,3+\la,3+\la^2\}$ representing the third ideal
class.

Multiplications by $\la_1$ and $\la_2$ generate the following
three weakly algebraically isomorphic actions which are not
algebraically isomorphic by Theorem \ref{mytheor} even up to a
time change, and therefore not measurably isomorphic:

\[
A=\left(
	\begin{smallmatrix} \hphantom{-}0 & \hphantom{-}1&
\hphantom{-}0\\ \hphantom{-}0 & \hphantom{-}0 &
\hphantom{-}1\\ -1 & \hphantom{-}8 & \hphantom{-}2
\end{smallmatrix}\right) \quad\text{and} \quad
B=\left(	\begin{smallmatrix} \hphantom{-}2
& \hphantom{-}1 & \hphantom{-}0\\ \hphantom{-}0 &
\hphantom{-}2 & \hphantom{-}1\\ -1 &
\hphantom{-}8 & \hphantom{-}4
\end{smallmatrix}\right);
\]
\[
A'=\left(
	\begin{smallmatrix} -1 &
\hphantom{-}2 &
\hphantom{-}0\\-1 & \hphantom{-}1 & \hphantom{-}1\\
-6 & \hphantom{-}9 & \hphantom{-}2
\end{smallmatrix}\right) \quad\text{and} \quad
B'=\left(	\begin{smallmatrix} \hphantom{-}1 & \hphantom{-}2 &
\hphantom{-}0\\-1 & \hphantom{-}3 & \hphantom{-}1\\
-6 & \hphantom{-}9 & \hphantom{-}4
\end{smallmatrix}\right);
\]
\[
A''=\left(
	\begin{smallmatrix}-3 & \hphantom{-}4 & \hphantom{-}0\\
-3 & \hphantom{-}3 & \hphantom{-}1\\
-10 & \hphantom{-}11 & \hphantom{-}2
\end{smallmatrix}\right) \quad\text{and} \quad
B''=\left(	\begin{smallmatrix} -1 & \hphantom{-}4 &
\hphantom{-}0\\ -3 & \hphantom{-}5 & \hphantom{-}1\\
-10 & \hphantom{-}11 & \hphantom{-}4
\end{smallmatrix}\right).
\]

Each action has $2$ fixed point in $\T^3$, $(0,0,0)$ and
$(\frac 1{2},\frac 1{2},\frac 1{2})$. Hence all affine
centralizers are isomorphic as abstract groups to $\Z^2\times
\Z/2\Z\times
\Z/2\Z$.

\medskip

\noindent{\bf Example 3c} Finally we give an example of two nonisomorphic
maximal  Cartan actions which come from  the vector of fundamental 
units $\bar\la=(\la_1,\la_2)$ in a totally real cubic field $K$ such 
that $\Z(\la_1,\la_2)\neq \mathcal O_K$.
Thus the whole group of units  does not generate the ring $\mathcal O_K$.
Both actions
  $\al_{\bar\la}^{\min}$ and  $\al_{\bar\la}^{\max}$  of the group $\Z^2$
are maximal Cartan actions by Lemma \ref{l:rings-units}. However by 
Corollary \ref{cor-cyclic} the former is cyclic
and the latter is not
  and hence they are not measurably isomorphic up to a time change by
Corollary  \ref{cor:min-max}.

For a specific example we pick the totally real cubic field 
$K=\Q(\al)$ with class number $1$ discriminant $1304$ given by the 
polynomial $x^3-x^2-11x-1$. For this filed we have $[\mathcal 
O_K:\Z(\al)]=2$. Generators in $\mathcal O_K$ can be taken to be 
$\{1,\al,\beta=\frac{\al^2+1}2\}$. Fundamental units are 
$\la_1=-\al,\;\la_2=-5+14\al+10\beta=14\al+5\al^2\in\Z[\al]$. Thus 
the whole group of units lies in $\Z[\la]$. To construct the 
generators for two non--isomorphic action  $\al_{\bar\la}^{\min}$ and 
$\al_{\bar\la}^{\max}$  we write multiplications by $\la_1$ and 
$\la_2$ in bases $\{1,\al,\al^2\}$ and $\{1,\al,\beta\}$, 
correspondingly. The resulting matrices are:
\[
A=\left(\begin{smallmatrix}\hphantom{-}0 & -1 & \hphantom{-}0\\
\hphantom{-}1 & \hphantom{-}0 & -1 \\
\hphantom{-}1 & \hphantom{-}11 & \hphantom{-}1
\end{smallmatrix}\right)\quad
B=\left(\begin{smallmatrix} 0 & 14 & 5\\ 5 & 55 & 19\\
19 & 214 & 74
\end{smallmatrix}\right),
\]
\[
A'=\left(\begin{smallmatrix}\hphantom{-}0 & -1 & \hphantom{-}0\\
\hphantom{-}1 & \hphantom{-}0 & -2 \\
\hphantom{-}0 & -6 & -1
\end{smallmatrix}\right)\quad
B=\left(\begin{smallmatrix} -5 & \hphantom{-}14 & \hphantom{-}10\\ 
-14 & \hphantom{-}55 & \hphantom{-}38\\
-30 & \hphantom{-}114 & \hphantom{-}79
\end{smallmatrix}\right).
\]
  The first action has only one fixed point, the origin; the second 
has four fixed points  $(0,0,0)$,
$(\frac 1{2},\frac 1{2},\frac 1{2})$, $(\frac 1{2},\frac
1{2},0)$, and $(0,0,\frac 1{2})$.  Thus we have an example of two 
maximal Cartan actions of $\Z^2$ which have nonisomorphic affine and 
hence measurable centralizers.

\end{document}